\newcommand{\ba}{\begin{array}}
\newcommand{\ea}{\end{array}}
\newcommand{\bc}{\begin{center}}
\newcommand{\ec}{\end{center}}

\newcommand{\beqn}[1]{\begin{equation}\label{#1}}
\newcommand{\eeqn}{\end{equation}}
\newcommand{\be}{\begin{equation}}
\newcommand{\ee}{\end{equation}}

\newcommand{\beqnn}{\begin{eqnarray}}
\newcommand{\eeqnn}{\end{eqnarray}}

\newcommand{\col}{{\rm col}}

\documentclass{ifacconf}
\usepackage{graphicx}      
\usepackage{natbib}        
\usepackage{romannum}
\usepackage{amsmath}
\usepackage{mathtools}
\usepackage{xcolor}
\usepackage{algorithm} 
\usepackage{booktabs}
\usepackage{algpseudocode} 
\DeclarePairedDelimiter{\norm}{\lVert}{\rVert}
\usepackage{amsfonts}
\usepackage{bm}
\usepackage{fancyhdr}
\begin{document}
\begin{frontmatter}
\title{Integrated motion control and energy management of series hybrid electric vehicles: A multi-objective MPC approach}%\thanksref{footnoteinfo}} 
% \thanks[footnoteinfo]{Sponsor and financial support acknowledgment
% goes here. Paper titles should be written in uppercase and lowercase
% letters, not all uppercase.}

\author[First]{Henglai Wei}
\author[Third]{Guangyuan Li}
\author[Second]{Yang Lu}
\author[Third]{Hui Zhang}
% \author[First]{Yang Shi} 
\address[First]{Centre of Excellence for Testing \& Research of Autonomous Vehicles NTU (CETRAN), Nanyang Technological University, 639798, Singapore (Email: henglai.wei@ntu.edu.sg).}
\address[Third]{School of Transportation Science and Engineering, Beihang University, Beijing, China, 100191 (Email: huizhang285@gmail.com).}
\address[Second]{China North Vehicle Research Institute, Beijing, China, 100072.}

\begin{abstract}         
This paper considers the integrated motion control and energy management problems of the series hybrid electric vehicles (SHEV) with constraints. We propose a multi-objective model predictive control (MOMPC)-based energy management approach, which is embedded with the motion control to guarantee driving comfort. In addition, due to the slow response of the engine, it may cause excessive batter power when HEVs work in different conditions (e.g., uphill or sudden acceleration) with a certain request power; this implies the discharge current is too large. A battery current constraint is designed and incorporated into the MOMPC optimization problem and hence avoids the extra high charge-discharge current. This prevents potential safety hazards and extends the battery's life. Finally, the numerical experiments are performed to verify the proposed approach.
\end{abstract}
\end{frontmatter}
%=============================================================================
\section{Introduction}
In recent years, electric vehicles have been recognized as a critical step towards energy conservation and emission reduction because of their environmental protection, energy saving, low noise and many other technical advantages \cite{wu2015electric}. Hybrid electric vehicles (HEV) combine the advantages of electric vehicles and traditional fuel vehicles, which is regarded as one of the effective ways to save energy and reduce emissions. In particular, hybrid electric vehicles (HEV) embody the attributes of high performance and low emissions \cite{hannan2014hybrid}. For these HEVs, a sound energy management strategy is the key to improving the vehicle's fuel economy. The energy management strategy can realize the reasonable distribution of multiple power sources (e.g., the internal combustion engine (ICE) and the battery), therefore, they can achieve the purpose of energy saving and emission reduction under different working conditions. According to the architecture and configuration of a hybrid electric powertrain, HEVs can be classified into different types, for example, series HEV (SHEV), parallel HEV (PHEV), and parallel-series HEV (PSHEV) \cite{ehsani2021state}. This paper focuses on the energy management problem of SHEVs. Though there has been some progress in energy management for HEVs \cite{sabri2016review}, developing a more efficient energy management strategy remains a challenge that deserves further research.

Being one of the main concerns of HEVs, plenty of results on energy management strategies have been reported in the literature, including the rule-based strategies \cite{hofman2007rule,banvait2009rule,trovao2013multi,peng2017rule}, optimization-based strategies \cite{ettihir2016optimization,chen2019series,hu2022multihorizon}, and learning-based strategies \cite{wu2018continuous,li2019deep,lian2020rule}. Over the past two decades, numerous optimization-based techniques have been devised to increase the fuel efficiency of Hybrid Electric Vehicles (HEVs). One such technique involves the use of dynamic programming (DP) \cite{koot2005energy, chen2014energy}, which guarantees global optimality. Nevertheless, the high computational complexity of DP-based energy strategies restricts their practical utility. Other commonly used approaches to enable energy awareness of HEVs in optimization-based strategies are to employ Pontryagin’s Minimum Principle (PMP) \cite{hou2014approximate,uebel2018optimal} or equivalent consumption minimization strategies (ECMS) \cite{musardo2005ecms}. Although these approaches are computationally efficient, it is hard to ensure optimality accuracy, especially when the system model complexity increases.  In fact, optimality accuracy and real-time implementation are required in most of the existing energy management strategies for practical HEVs.

Model predictive control (MPC) is advocated here because of its excellent ability to handle physical constraints and obtain high performance \cite{wei2021distributed,wei2022self}. Some interesting results on MPC-based energy management strategies have been developed for HEVs \cite{borhan2012mpc,di2013stochastic,wang2016model}. On the other hand, MPC has also been utilized to address the motion control problem of the single autonomous vehicle \cite{beal2013model} or cooperative autonomous vehicles \cite{zheng2017distributed}. In this paper, we propose to enhance driving comfort in energy management strategies by designing a motion control-related cost function. More specifically, we incorporate the motion control cost function into the overall cost function design of the energy management task. It is worth noting that the above-mentioned MPC-based energy management approaches only consider the decoupled powertrain level dynamics. %Simultaneously guaranteeing driving comfort and improving the fuel efficiency of HEVs remains challenging. 

It is worth noting that the engine subsystem of a SHEV has special characteristics. When the HEVs work in different conditions (e.g., the acceleration and starting working conditions), the slow response of the engine generally leads to the battery's power output is too high. The excessive discharge current affects the battery's efficiency and life \cite{han2019review}. Although the approaches in \cite{chen2019series,zheng2020integrated} address motion control and energy management simultaneously, they do not consider the problem of excessive battery current, which, however, is the key to guaranteeing the battery's safety and life. To address this challenging issue, we design a box constraint to bound the current state of the battery, which is distinct from the existing co-optimization approaches \cite{chen2019series,zheng2020integrated}. 

The main contributions of our work are summarized as follows: 1) A multi-objective MPC (MOMPC) approach is proposed to solve the integrated motion control and energy management issues in constrained HEVs. By addressing these challenges together, the proposed approach can optimize the joint optimization problem, leading to enhanced driving comfort and fuel efficiency. 2) A current constraint is devised for the battery to prevent excessive charge-discharge current during different HEV operating conditions, including acceleration. This current constraint is a novel feature of our approach, distinguishing it from existing solutions, which often overlook this aspect. 3) We perform simulations of motion control and energy management problems on a Sample HEV to validate the effectiveness of our approach. Overall, The proposed MOMPC approach contributes to the optimization of practical HEVs, leading to enhanced fuel energy efficiency, driving comfort, and reduced environmental impact.

% The paper is organized as follows: Section \ref{whl14a-sec:2} formulates the motion control problem and energy management problem of a SHEV with constraints. In Section \ref{whl14a-sec:3}, we propose the MOMPC approach. Section \ref{whl14a-sec:4} presents the simulation results. Finally, Section \ref{whl14a-sec:5} concludes the paper.

\textbf{Notation}: For any vector ${x} \in \mathbb{R}^n$, $\norm{x}_{P}^2$ denotes the weighted norm ${x}^{{T}}{Px}$. $[{x}_1^\text{T},\dots,{x}_n^\text{T}]^\text{T}$ is written as $\col({x}_1,\dots,{x}_n)$. ${x}(t)$ denotes the state ${x}$ at time $t$, and ${x}(t+k;t)$ denotes the predicted state at some future time $t+k$ determined at time $t$.

\section{Problem Formulation}
\label{whl14a-sec:2}
%The superscript `${T}$' denotes the transposition.

\subsection{Vehicle longitudinal dynamics}
The vehicle longitudinal dynamics are described by 
\begin{equation}\label{whl14a-eq:1}
\begin{cases}
\dot{s}=v,\\
\dot{v}=\dfrac{1}{m}\big(\dfrac{\eta F_{d}}{R_w}-C_Av^2-mgf\cos \theta-mg\sin\theta\big),
\end{cases}
\end{equation}
where $s$ and $v$ are, respectively, the position and the velocity along the longitudinal axis; $\eta$ is the mechanical efficiency of the driveline, $R_w$ is the tire radius, $C_A$ is the aerodynamic drag coefficient, $m$ is the vehicle mass, $g$ is the gravity constant, $f$ is the rolling resistance, $\theta$ is the road slope, and $F_{{d}}$ is the desired driving torque. Note that the vehicle is enforced to satisfy the state and control input constraints
\begin{equation}
0\leq v\leq v^{\max}, \ F^{\min}\leq F_d\leq F^{\max}.
\end{equation}

Next, the vehicle longitudinal dynamics \eqref{whl14a-eq:1} is discretized first by applying the Euler forward method for the MOMPC formulation, i.e.,
\begin{equation}
x(t+1)=Ax(t)+Bu(t),
\end{equation}
where $x(t)=[s(t),v(t)]^T$, $A=[1,\delta t;0,1]$, $B=[0;1]$, and $\delta t$ is the sampling period.

\subsection{Engine and battery model}
The total wheel power $P_v$ is given as follows
\begin{equation}
P_v=F_dv. 
\end{equation}
In the following, the main components (i.e., the engine and the battery) of a power-split SHEV are modeled for the energy management and motion control purpose. For SHEV, the requested power $P_{r}$ is provided by the engine and the battery, i.e., 
\begin{equation}
P_v=P_r=\eta_m(P_e+P_b),
\end{equation}
where $P_e$ is the engine power, $P_b$ is the batter power, and $\eta_m$ is the motor efficiency.

\begin{figure}[!ht]
\includegraphics[width=1\columnwidth]{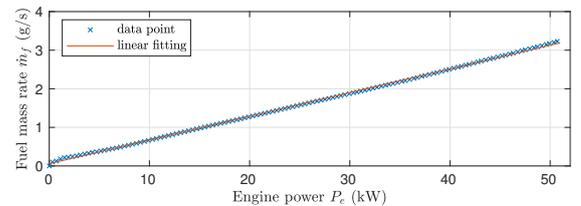}
\centering
\caption{The relationship between fuel mass rate $\dot{m}_f$ and engine power $P_e$.}
\label{whl14a-fig:2}
\end{figure}

\emph{Engine model:} The engine power is determined by
\begin{equation}
P_e=w_e\tau_e,
\end{equation}
where $w_e$ is the engine speed and $\tau_e$ is the engine torque. In order to improve fuel efficiency, the engine needs to be operated at the most efficient power-speed operating points. Note that the efficient operation trajectory is pre-defined. That is, once the engine power is given, we can find the corresponding torque and speed of the engine. As indicated in the empirical map of the engine, the relationship between the fuel flow rate $\dot{m}_f$, the engine speed and torque can be represented by a nonlinear function $\psi$, that is
\begin{equation}
\frac{d}{dt}{m}_f=\psi(w_e,\tau_e).
\end{equation}
Moreover, the fuel mass dynamics can be filtered by
\begin{equation}\label{whl14a-eq:7}
\frac{d}{dt}{m}_f=\alpha P_e+\beta,
\end{equation}
in which the engine power coefficient $\alpha$ and the idle fuel mass rate $\beta$ are obtained via linear regression method. As shown in Fig.~\ref{whl14a-fig:2}, the relationship between the fuel mass rate and the engine power $P_e$ is approximately linear. The fuel mass dynamics take the following discrete-time form
\begin{equation}
\Delta{m}_f=\delta t(\alpha P_e+\beta).
\end{equation}

Meanwhile, the engine torque and speed should satisfy the physical constraints
\begin{equation}
w_e^{\min} \leq w_e\leq w_e^{\max},\ \tau_e^{\min} \leq \tau_e\leq \tau_e^{\max}.
\end{equation}

\emph{Battery model:} The battery's state of charge (SOC) dynamics are governed by
\begin{equation}
\frac{d}{dt}{SOC}=-I_b/Q_b,
\end{equation}
in which $I_b$ and $Q_b$ denote the current and capacity of the battery, respectively. The battery closed circuit voltage is estimated by $V_b=V_{oc}-I_bR_b$, where $V_{oc}$ denotes the open circuit voltage and $R_b$ denotes the battery resistance. The corresponding power $P_b$ provided by the battery is described by $P_b=I_bV_b$. Then, we get
\begin{equation}
I_b=\frac{V_{oc}-\sqrt{V_{oc}^2-4R_bP_b}}{2R_b}.
\end{equation}

Then, we obtain 
\begin{equation}\label{whl14a-eq:13}
\frac{d}{dt}SOC=-\frac{V_{oc}-\sqrt{V_{oc}^2-4R_bP_b}}{2R_bQ_b}.
\end{equation} 
Moreover, the battery SOC model \eqref{whl14a-eq:13} is discretized first by applying the Euler forward method for the MOMPC formulation, i.e.,
\begin{equation}
SOC(t+1)=SOC(t)-\frac{V_{oc}-\sqrt{V_{oc}^2-4R_bP_b}}{2R_bQ_b}\delta t.
\end{equation}

In order to improve the battery life, the battery's SOC and current are enforced to be operated in reasonable ranges 
\begin{equation}
SOC^{\min}\leq SOC\leq SOC^{\max},\ I_b^{\min}\leq I_b\leq I_b^{\max}.
\end{equation}

% Furthermore, the control-oriented vehicle dynamics and energy model can be characterized by
% \begin{equation}
% x
% \end{equation}

\subsection{Problem formulation}
This work aims to develop an effective MOMPC strategy enabling the codesign of energy management and motion control for the SHEV with constraints. In particular, discharging under high current may lead to accelerated degradation and safety problems of the battery \cite{han2019review}. The proposed strategy can minimize the fuel consumption by finding a suitable energy power split and enhance driving comfort while avoiding the excessive battery current of the HEV under different working conditions.

\section{Main results}
\label{whl14a-sec:3}
\begin{figure}[ht]
\includegraphics[width=0.8\columnwidth]{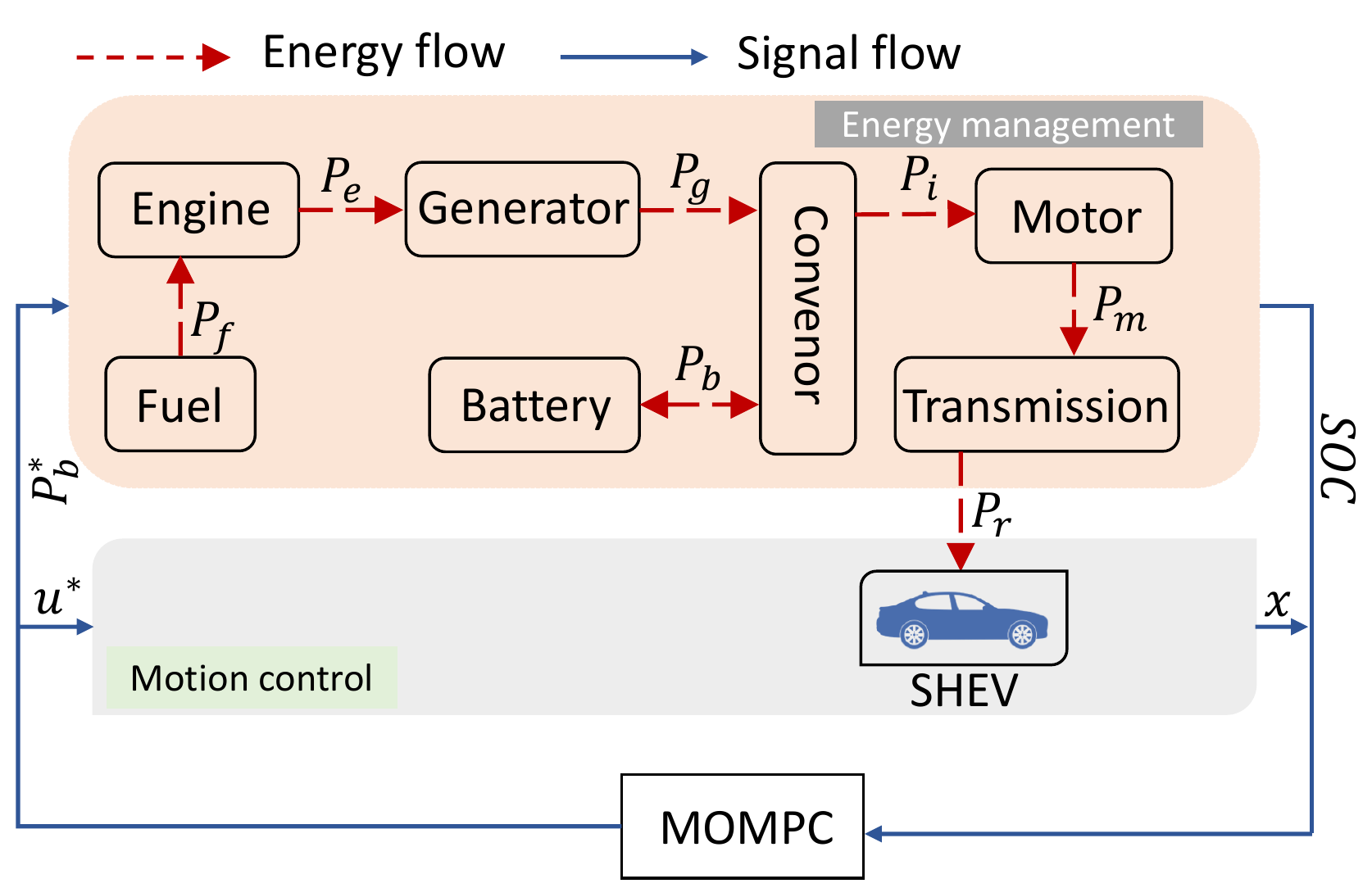}
\centering
\caption{Block diagram of the proposed MOMPC for the energy management and motion control of the SHEV.}
\label{whl14a-fig:1}
\end{figure}

In order to simultaneously overcome these challenges of the SHEV, a MOMPC strategy that simultaneously considers the energy management and motion control problems is proposed in this section. Fig.~\ref{whl14a-fig:1} shows the architecture of the proposed MOMPC strategy for the energy management and motion control of a SHEV. 

\graphicspath{{./Fig/}}

\subsection{Objective function}
% Different control strategies for the motion control problem of the conventional vehicle has been well studied in the literature {\color{red}[X]}. However, these solutions may not be suitable for the SHEV. 
Before formulating the MOMPC optimization problem, we design the objective functions concerning the motion control and energy consumption.

\emph{1) Cost function for motion control:} Given the reference distance and speed $x_d$, the cost function for motion control for the SHEV at time $t$ is expressed as the following
\begin{equation}
{J}_{m}=\sum_{k=0}^{N}\|{x(k;t)-x_d(k;t)}\|_{Q}^2+\sum_{k=0}^{N-1}\|u(k;t)\|_{S}^2
\end{equation}
where $N$ denotes the prediction horizon, $Q=Q^T \succeq 0$, $S=S^T \succ 0$, ${x}_d(k;t)$ is the desired state sequence, and the motion control cost function is assumed to be continuous. 

\emph{2) Cost function for fuel consumption:} For the SHEV, the cost function for the fuel consumption at time $t$ is expressed as the following
\begin{equation}
{J}_{f}=\sum_{k=0}^{N}\|\Delta m_f(k;t)\|_{R}^2
\end{equation}
where $\Delta m_f(k;t)$ denotes the predicted energy consumption at time $t+s$ and $R=R^T \succeq 0$.

\emph{3) Cost function for battery management:} The cost function for the battery management at time $t$ is given by
\begin{equation}
{J}_{b}=\sum_{k=0}^{N}\|SOC(k;t)-SOC_r\|_{P}^2
\end{equation}
where $SOC(k;t)$ represents the predicted SOC state in the future time $t+s$, $SOC_r$ is the reference state for the battery SOC, and $P=P^T \succeq 0$.

Note that the conventional SHEV's energy management strategy aims to minimize fuel consumption while enforcing the practical SOC state close to the reference SOC state \cite{borhan2012mpc,wang2016model}.
% \begin{rem}
% Different from the consistency constraint in \cite{zhan2018distributed,zheng2016distributed}, the local stage cost term $\norm{x_i(s;t_{k_i}^i)-{x}_i^b(s+H^{i*}(t^i_{k_i-1});t_{k_i-1}^{i})}_{Q_i'}^2$ enforces a degree of consistency between what subsystem $i$ plans to do and what neighbors believe subsystem $i$ will do. In this way, the constraint is transformed into a soft constraint, and the feasible solution is easy to construct.      
% \end{rem}
\subsection{MOMPC optimization}
In what follows, the MOMPC optimization problem for the SHEV is defined. 
\subsubsection{Pareto MOMPC optimization}
The Pareto MOMPC optimization problem $\mathcal{P}_1$ for the SHEV at time instant $t$ is formulated in the following
\begin{equation}
\begin{aligned}
\min_{\textbf{u}(t)}\ &J_p=[{J}_{1},{J}_{2}]^T\\
\text{s.t.}\ &x(k+1;t)=Ax(k;t)+Bu(k;t),\\
             &P_r(k;t) = \eta_m(P_e(k;t)+P_b(k;t)),\\
            &SOC(k+1;t)=SOC(k;t)-\frac{I_b(k;t)}{Q_b}\delta t,\\
            &I_b(k;t)=\frac{V_{oc}-\sqrt{V_{oc}^2-4R_bP_b(k;t)}}{2R_b},\\
            &SOC(0;t)=SOC(t),\ {x}(0;t)={x}(t), \\
            &SOC^{\min}\leq SOC(k;t)\leq SOC^{\max},\\
            &I_b^{\min}\leq I_b(k;t)\leq I_b^{\max},\\
            &{\Delta m_f=(\alpha P_e(k;t)+\beta)\delta t,}\\
            &w_e^{\min} \leq w_e(k;t)\leq w_e^{\max},\tau_e^{\min} \leq \tau_e(k;t)\leq \tau_e^{\max},
      % &x_i(s;t_{k_i}^i)\in \mathbb{X}_i\label{whl22-eq:5d}
\end{aligned}
\end{equation}
in which $J_1=J_m$, $J_2=J_f+J_b$, $\textbf{u}(t)=[\bm{u}(t),\bm{P}_b(t)]^T$ denotes the control input sequence to be optimized, with $\bm{u}(t)=\col(u(0;t),\dots,u(N-1;t))$ and $\bm{P}_b(t)=\col(P_b(0;t),\dots,P_b(N-1;t))$. Here, the objective function $J_p=[J_1,J_2]^T$ is a vector-valued function, which consists of two different objective functions. 

It is impossible to simultaneously optimize two objectives. In the following, the Pareto optimality \cite{deb2014multi} is introduced. A solution $\textbf{u}^*$ is said to be Pareto another feasible solution $\textbf{u}$ if the conditions $J_i(\textbf{u}^*)\leq J_i(\textbf{u}),\forall i\in[1,2]$ and $J_i(\textbf{u}^*)< J_i(\textbf{u}),\exists i\in[1,2]$ hold. It is worth noting that there generally exists more than one Pareto optimal solution. Specially, these objective vectors under the Pareto optimal solutions are represented by a Pareto frontier. These solutions and their corresponding objective function values form the Pareto set, which helps to find the most compromised Pareto optimal solution. The evolutionary methods \cite{deb2014evolutionary} can be employed to determine the whole Pareto frontier. However, it is not necessary to determine the whole Pareto frontier for the energy management and motion control problems of the SHEV. The weighted-sum MOMPC is adopted to find one preferred solution from the Pareto frontier in this work.

\subsubsection{Weighted-sum MOMPC optimization}
The weighted-sum MOMPC method combines different objective functions into one single objective function by assigning configurable weight to each objective. 
The weighted-sum MOMPC optimization problem $\mathcal{P}_2$ for the SHEV at time instant $t$ is formulated in the following
\begin{subequations}\label{whl14a-eq:22}
\begin{alignat}{2}
\min_{\textbf{u}(t)}\ & J=\alpha_1{J}_{m}+\alpha_2{J}_{f}+\alpha_3{J}_{b}\\
\text{s.t.}\ &x(k+1;t)=Ax(k;t)+Bu(k;t),\\
             &P_r(k;t) = \eta_m(P_e(k;t)+P_b(k;t)),\\
            &SOC(k+1;t)=SOC(k;t)-\frac{I_b(k;t)}{Q_b}\delta t,\\
            &I_b(k;t)=\frac{V_{oc}-\sqrt{V_{oc}^2-4R_bP_b(k;t)}}{2R_b},\label{whl14a-eq:20h}\\
            &SOC(0;t)=SOC(t),\ x(0;t)=x(t), \\
            &SOC^{\min}\leq SOC(k;t)\leq SOC^{\max},\\
            &I_b^{\min}\leq I_b(k;t)\leq I_b^{\max},\\
            &{u}(k;t)\in \mathcal{U},\ x(k;t)\in \mathcal{X},\\
            &{\Delta m_f=(\alpha P_e(k;t)+\beta)\delta t,}\\
            &w_e^{\min} \leq w_e(k;t)\leq w_e^{\max},\ \tau_e^{\min} \leq \tau_e(k;t)\leq \tau_e^{\max},
\end{alignat}
\end{subequations}
where $\alpha_i$ denotes the cost weight with $\sum_{i=1}^3\alpha_i=1$. $\mathbf{u}^*(t)$ is generated by calculating the optimization problem $\mathcal{P}_2$ at time $t$. $\tilde{\mathbf{u}}(t)$ is a feasible control input sequence. Note that the weight parameter $\alpha_i$ represents the preference of the decision maker, which balances the trade-off between the motion control performance and fuel energy consumption for the SHEV. The weight parameter can be chosen by a logistic function as in \cite{shen2018path}. 

Eq. \eqref{whl14a-eq:20h} indicates that the excessive battery current of the SHEV under special working conditions is handled by incorporating the hard constraint into the MOMPC optimization problem. Furthermore, the speed profile prediction based on historical traffic data can be exploited to improve fuel consumption efficiency and avoid excessive battery current \cite{xiang2017energy}. Suppose the acceleration of the HEV can be accurately predicted. In that case, we can turn on the engine in advance to avoid excessive battery discharge current, which achieves the purpose of protecting the battery. 

\subsection{Weighted-sum MOMPC algorithm}
The proposed weighted-sum MOMPC algorithm for the SHEV is specified as follows. 
\begin{algorithm}
\caption{Weighted-sum MOMPC Algorithm}
\label{alg:Framwork}   
\begin{algorithmic}[1] 
\State \textbf{Initialization}: For the SHEV, give the initial states $x(t)$, $SOC(t)$, the initial feasible control $\tilde{\textbf{u}}(t)$ and other design parameters. Set $t=0$. 
\State Sample system state $x(t)$;\label{step2}
\State Solve the optimization problem $\mathcal{P}_2$ \eqref{whl14a-eq:22} and generate the optimal control $\mathbf{u}^{*}(t)$; 
\State Distribute the required power $P_r(t)$ between the battery $P_b(t)$ and engine $P_e(t)$;
\State Apply the control input $u^*(t)$ to the SHEV;
\State $t=t+1$ and go to Step $2$; if $t=T_{\text{sim}}$, stop.
\end{algorithmic}  
\end{algorithm}  

\section{Numerical simulations}
\label{whl14a-sec:4}
This section evaluates the effectiveness of the proposed approach using numerical simulations over the Urban Dynamometer Driving Schedule (UUDS) driving cycle. The MOMPC optimization problem $\mathcal{P}_2$ in \eqref{whl14a-eq:22} is solved numerically using Yalmip, which exploits IPOPT algorithm for the numerical optimization.
The vehicle model parameters are summarized in Table~\ref{whl14a-tab:1}.
\begin{table}[!ht]
\begin{center}
\caption{SHEV model parameters.}\label{whl14a-tab:1}
\begin{tabular}{lll}
\toprule
Symbol  &  Value&  Description \\
\hline
$m$  & 1405kg & vehicle mass \\
$\eta$   & 0.96 &driveline efficiency\\
$R_w$  & 0.3050m &tire radius\\
$C_A$  & $0.5063$ &aerodynamic drag coefficient\\
$f$  & 0.01 &rolling resistance\\
$\eta_m$  & 0.96 &motor efficiency\\
$V_{oc}$  & $220.64$V &battery open circuit voltage\\
$R_b$  & $0.3757\Omega$ &battery resistance\\
$Q_b$  & $23.4$Ah &battery capacity\\
$SOC^{\min,\max}$  &$0.3$, $0.8$  &battery SOC limits\\
$I_b^{\min,\max}$  &$-90$, $90$A  &battery current limits\\
$w_e^{\min,\max}$& $0$, $105$rad/s& engine speed limits\\
$\tau_e^{\min,\max}$  & $0$, $112$N/m &engine torque limits\\
\bottomrule
\end{tabular}
\end{center}
\end{table}

The engine power coefficient and the idle fuel mass rate are identified as $\alpha=0.0614$ and $\beta=0.0583$, respectively. The sampling time is chosen as $\delta t=1s$, and the prediction horizon is $N=10s$. The reference SOC is set as $SOC_r=0.5$. The weighting matrices for the cost functions are selected as $Q=1$, $S=1$, $R=5$ and $P=300$, and the cost weights are $\alpha_1=0.33$, $\alpha_2=0.33$, $\alpha_3=0.33$. The initial states of the SHEV are $x(0)=[0, 0]^{T}$, $SOC(0)=0.66$.

\begin{figure}[!ht]
\includegraphics[width=1\columnwidth]{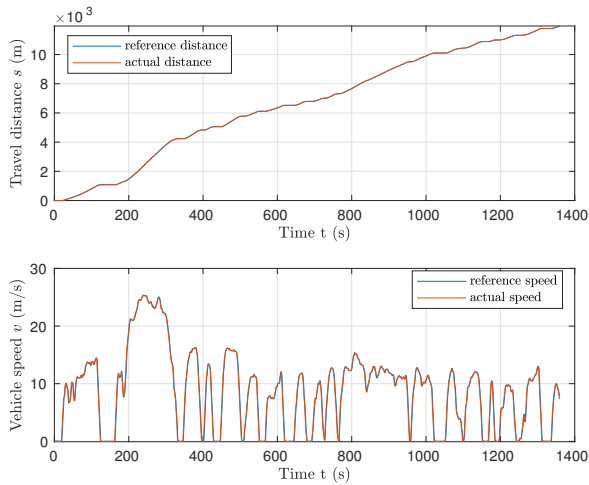}
\centering
\caption{State trajectories of the SHEV under the proposed MOMPC algorithm. ({Top}): The traveled distance of the SHEV. ({Bottom}): The SHEV velocity.}
\label{whl14a-fig:3}
\end{figure}

\begin{figure}[!ht]
\includegraphics[width=1\columnwidth]{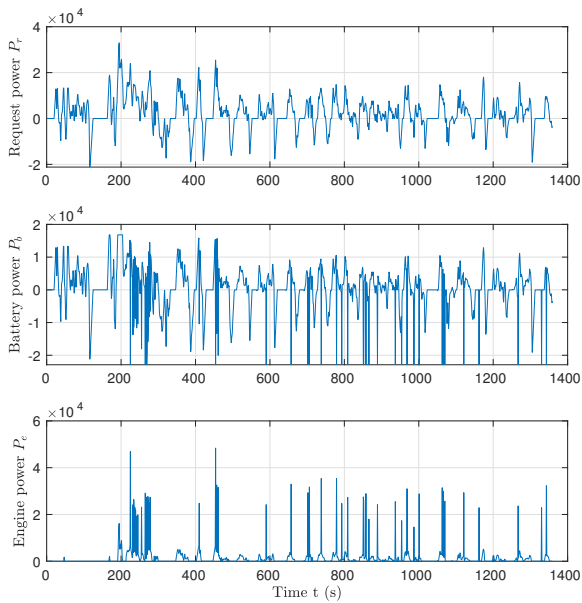}
\centering
\caption{UUDS cycle simulation. (Top): The requested power $P_r$ [W]. (Middle): The power is delivered by the battery $P_b$ [W]. (Bottom): The power supplied by the engine $P_e$ [W].}
\label{whl14a-fig:4}
\end{figure}

\begin{figure}[!ht]
\includegraphics[width=1\columnwidth]{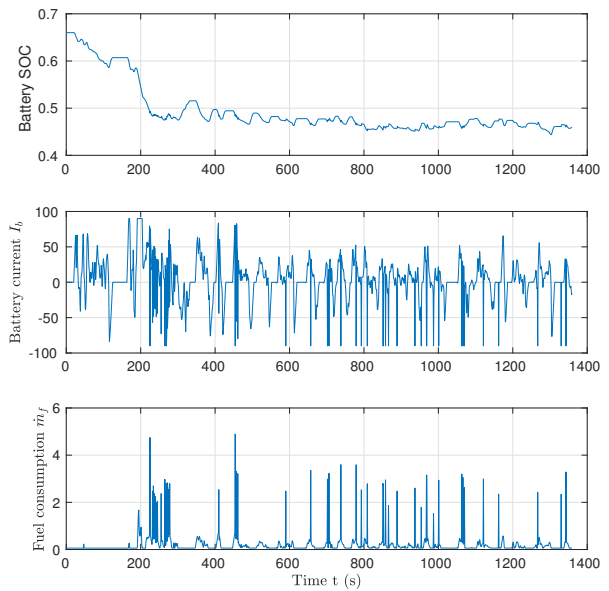}
\centering
\caption{UUDS cycle simulation. (Top): The battery SOC. (Middle): The battery current $I_b$ [A]. (Bottom): The fuel consumption rate $\dot{m}_f$ [g/s].}
\label{whl14a-fig:5}
\end{figure}

The trajectory tracking and energy management simulation results are presented in Fig.~\ref{whl14a-fig:3} - Fig.~\ref{whl14a-fig:5}. The traveled distance and velocity of the SHEV are depicted in Fig.~\ref{whl14a-fig:3}. It can be observed that the vehicle can successfully track the reference trajectory. Fig.~\ref{whl14a-fig:4} illustrates the power split of the SHEV, including the required power, battery power and engine power, implying that the required power can be well split between the battery and engine. Fig.~\ref{whl14a-fig:5} shows the profiles of the battery SOC, the battery current and the fuel consumption rate over the UUDS driving cycles. As we can see, 1) the SOC is maintained within reasonable limits; 2) the battery current is always within the corresponding permitted ranges. It's worth noting that the negative current shown in Fig.~\ref{whl14a-fig:5} indicates that the battery is being charged. Consequently, the corresponding battery power, $P_b$, is also negative as shown in Fig.~\ref{whl14a-fig:4}.

% \graphicspath{{./Fig/}}
% \begin{figure*}[ht]
% \includegraphics[width=2\columnwidth]{Fig5_disturbances.eps}
% \centering
% \caption{Disturbances trajectories of three subsystems.}
% \label{whl22-fig:5}
% \end{figure*}

% \begin{table}[!ht]
% \begin{center}
% \caption{Performance comparison.}\label{whl22-tab:1}
% \begin{tabular}{ccc}\hline\hline
% Method  &  Average sampling time&  Average cost  \\\hline
% Periodic  & {0.3000} & 1.4702 \\
% self-triggered   & \bf{0.9065} &1.5600\\ \hline
% \end{tabular}
% \end{center}
% \end{table}
%%%%%%%%%%%%%%%%%%%%%%%%%%%%%%%%%%%%%%%%%%%%%%%%%%%%%%%%%%%%%%%%%%%%%%%%%%%%%%
        %%%%%%%%%%%%%%%%%%%%%%% Section V %%%%%%%%%%%%%%%%%%%%%%%%%%%%%
%%%%%%%%%%%%%%%%%%%%%%%%%%%%%%%%%%%%%%%%%%%%%%%%%%%%%%%%%%%%%%%%%%%%%%%%%%%%%%
\section{Conclusion}
\label{whl14a-sec:5}
In this paper, we have proposed a MOMPC approach for the integrated motion control and energy management problem of the SHEV with constraints. The proposed approach can simultaneously enhance driving comfort and improve fuel consumption efficiency. Generally, the engine responds more slowly than the battery when the SHEV demands more power under specific working conditions (e.g., uphill or sudden acceleration). In this case, the sharply increased power demand may cause an excessive discharge current. A battery current constraint was designed and incorporated into the MOMPC optimization problem, avoiding the extra high charge-discharge current. The proposed MOMPC approach guaranteed the battery's safety and extended the battery's life. 
The simulation results verified the effectiveness of the proposed approach. Future research will consider trajectory prediction in the energy management task.%, based on which the power demand can be predicted, and different power sources can be well distributed. 

\bibliography{ifacconf}  
\end{document}